\documentclass[a4paper,12pt]{article}
\usepackage[T1]{fontenc}
\usepackage[latin1]{inputenc}
\usepackage{amsthm,amsmath}
\usepackage{amssymb}
\usepackage[all]{xy}
\usepackage{enumerate}
\usepackage{a4wide}
\usepackage{hyperref}
\usepackage{latexsym}
\usepackage{enumerate}
\usepackage[mathscr]{eucal}

\newtheorem{theorem}{Theorem}[section]
\newtheorem{proposition}[theorem]{Proposition}
\newtheorem{corollary}[theorem]{Corollary}

\newtheorem*{theorem*}{Theorem}
\newtheorem*{proposition*}{Proposition}
\newtheorem*{corollary*}{Corollary}
\newtheorem*{lemma*}{Lemma}
\theoremstyle{definition}
\newtheorem{definition}[theorem]{Definition}

\newtheorem{example}[theorem]{Example}
\newtheorem{remark}[theorem]{Remark}
\newtheorem*{remark*}{Remark}

\newtheorem*{definition*}{Definition}
\numberwithin{equation}{section}

\newcommand{\cat}[1]{\mathcal{#1}}
\newcommand{\coring}[1]{\mathfrak{#1}}
\newcommand{\tensor}[1]{\otimes_{#1}}

\newcommand{\tensork}[2]{#1\otimes_{K} #2}
\newcommand{\rcomod}[1]{ \mathsf{Comod}_{#1}}
\newcommand{\rmod}[1]{\mathsf{Mod}_{#1}}
\newcommand{\lmod}[1]{{}_{#1}\mathsf{Mod}}

\newcommand{\cotensor}[1]{\square_{#1}}
\newcommand{\lcomod}[1]{{}_{#1}\mathsf{Comod}}
\newcommand{\homcom}[3]{\mathrm{Hom}^{#1}(#2,#3)}
\renewcommand{\hom}[3]{\mathrm{Hom}_{#1}\left( #2 \, , \, #3\right)}

\newcommand{\lend}[2]{\mathrm{End}({}_{#1}#2)}

\newcommand{\dkmod}[2]{{}_{#1}\mathcal{M}(\h)^{#2}}
\newcommand{\T}{\mathcal{T}}
\newcommand{\h}{\mathbf{H}}

\newcommand{\fk}[1]{\mathfrak{#1}}
\newcommand{\Sf}[1]{\mathsf{#1}}
\newcommand{\Bf}[1]{\mathbf{#1}}
\newcommand{\esc}[2]{\langle #1,#2 \rangle}

\newcommand{\Rat}{\mathrm{Rat}}

\newcommand{\lr}[1]{\left(\underset{}{} #1 \right)}
\newcommand{\Ker}[1]{\mathsf{Ker}(#1)}
\newcommand{\TT}{\mathscr{T}}

\begin{document}
\title{Corings with exact rational functors and injective objects \footnote{Investigaci\'{o}n realizada en el proyecto
MTM2004-01406 <<M\'{e}todos algebraicos en Geometr\'{\i}a no conmutativa>>
financiado por DGICYT y FEDER, y el proyecto P06-FQM-01889 <<Modelos
algebraicos aplicados a la F\'{\i}sica, Geometr\'{\i}a no conmutativa y
Computaci\'{o}n>> financiado por la Junta de Andaluc\'{\i}a}
\\ \vspace{0.15cm} \normalsize{\emph{Gewidmet an Professor Robert Wisbauer}}}
\author{L. El Kaoutit \\
\normalsize Departamento de \'{A}lgebra \\ \normalsize Facultad de Educaci\'{o}n y Humanidades \\
\normalsize  Universidad de Granada \\ \normalsize El Greco N$^o$
10, E-51002 Ceuta, Espa\~{n}a \\ \normalsize
e-mail: \textsf{kaoutit@ugr.es} \and J. G\'omez-Torrecillas \\
\normalsize Departamento de \'{A}lgebra \\ \normalsize Facultad de Ciencias \\
\normalsize Universidad
de Granada\\ \normalsize E18071 Granada, Espa\~{n}a \\
\normalsize e-mail: \textsf{gomezj@ugr.es} }

\date{ }
\maketitle

\begin{abstract}
We describe how some aspects of abstract localization on module
categories have applications to the study of injective comodules
over some special types of corings. We specialize the general
results to the case of Doi-Koppinen modules, generalizing previous
results in this setting.
\end{abstract}

\section*{Introduction}

The Wisbauer category $\sigma[M]$ subgenerated by a module $M$
\cite{Wisbauer:1991} is a flexible and useful tool when applied to
some at a first look unrelated situations. This has been the case of
the categories of comodules over corings, which, under suitable
conditions, become Wisbauer's categories
\cite{Abuhlail:2003a,Brzezinski/Wisbauer:2003,ElKaoutit/Gomez/Lobillo:2004c}.
On the other hand, as it was explained in \cite{Brzezinski:2002},
the categories of entwined modules and, henceforth, of Doi-Koppinen
modules, are instances of categories of comodules over certain
corings, which ultimately enlarges the field of influence of the
methods from Module Theory developed in \cite{Wisbauer:1991}. The
present paper has been deliberately written from this point of view,
although with a necessarily different style. To illustrate how
abstract results on modules may successfully applied to more
concrete situations, we have chosen a topic from the theory of
Doi-Koppinen modules with roots in the theory of graded rings and
modules, namely, the transfer of the injectivity from relative
modules (Doi-Koppinen, graded) to the underlying modules over the
ground ring (comodule algebra, graded algebra). This was studied at
the level of Doi-Koppinen modules in
\cite{Dascalescu/Nastasescu/Torrecillas:1999}, giving versions in
this framework of results on graded modules from
\cite{Dascalescu/alt:1996}. The methods developed in
\cite{Dascalescu/Nastasescu/Torrecillas:1999} rest on the exactness
of the rational functor for semiperfect coalgebras over fields
\cite{Gomez/Nastasescu:1995}, which allows the construction of a
suitable adjoint pair between the category of Doi-Koppinen modules
and the category of modules over the smash product \cite[Theorem
3.5]{Dascalescu/Nastasescu/Torrecillas:1999}. The pertinent
observation here, from the point of view of corings, is that one of
the functors in that adjoint pair is already a rational functor for
the coring associated to the comodule algebra \cite[Proposition
3.21]{Abuhlail:2003a}. Thus, a relevant ingredient in
\cite{Dascalescu/Nastasescu/Torrecillas:1999} is, under this
interpretation, the exactness of the trace functor defined by a
Wisbauer category of modules or, equivalently, the exactness of the
preradical associated to a closed subcategory of a category of
modules. Here, we make explicit the fact that the exactness of such
a preradical is equivalent to the property of being, up to an
equivalence of categories, the canonical functor of a localization
(Theorem \ref{TeoremaI}), and, henceforth, it has a right adjoint,
which is explicitly described. This right adjoint will preserve,
since it is a section functor, injective envelopes (Proposition
\ref{A-injectivoI}). We then deduce the general form of the transfer
of injective objects stated in
\cite{Dascalescu/Nastasescu/Torrecillas:1999}.

In the rest of this paper, we specialize the former general scheme
to corings with exact rational functors and, even more, to
Doi-Koppinen modules where the coacting coalgebra has an exact
rational functor.

The results of this paper should be not considered as completely
new. In fact, most part of them could be gathered, with suitable
adaptations (not always obvious), from other sources. Thus, our
text resembles a mini-survey. However, we believe that the reader
will not find elsewhere the statements made here, nor the
applications to the transfer of injectivity, since they do not
intend to be reproductions of previously published results. We
hope we have presented a study of some aspects of the theory of
corings and their comodules in a new light.

\medskip

\noindent \textbf{Notations and basic notions.} Throughout this
paper the word ring will refer to an associative unital algebra
over a commutative ring $K$. The category of all left modules over
a ring $R$ will be denoted by $\lmod{R}$, being $\rmod{R}$ the
notation for the category of all right $R$--modules. The notation
$X \in \cat{A}$ for a category $\cat{A}$ means that $X$ is an
object of $\cat{A}$, and the identity morphism attached to any
object $X$ will be denoted by the same character $X$.

Recall from \cite{Sweedler:1975} that an $A$--coring is a
three-tuple
$(\coring{C},\Delta_{\coring{C}},\varepsilon_{\coring{C}})$
consisting of an $A$--bimodule $\coring{C}$ and two homomorphisms
of $A$--bimodules (the comultiplication and the counity)
$$\xymatrix@C=50pt{\coring{C} \ar@{->}^-{\Delta_{\coring{C}}}[r] &
\coring{C}\tensor{A}\coring{C}},\quad \xymatrix@C=30pt{ \coring{C}
\ar@{->}^-{\varepsilon_{\coring{C}}}[r] & A}$$ such that
$(\Delta_{\coring{C}}\tensor{A}\coring{C}) \circ
\Delta_{\coring{C}} = (\coring{C}\tensor{A}\Delta_{\coring{C}})
\circ \Delta_{\coring{C}}$ and
$(\varepsilon_{\coring{C}}\tensor{A}\coring{C}) \circ
\Delta_{\coring{C}}=(\coring{C}\tensor{A}\varepsilon_{\coring{C}})
\circ \Delta_{\coring{C}}= \coring{C}$.

A right $\coring{C}$--comodule is a pair $(M,\rho_{M})$ consisting
of a right $A$--module $M$ and a right $A$--linear map $\rho_{M}:
M \rightarrow M\tensor{A}\coring{C}$, called right
$\coring{C}$--coaction, such that
$(M\tensor{A}\Delta_{\coring{C}}) \circ \rho_M =
(\rho_M\tensor{A}\coring{C}) \circ \rho_M$ and
$(M\tensor{A}\varepsilon_{\coring{C}}) \circ \rho_M=M$. A morphism
of right $\coring{C}$--comodules (or a right $\coring{C}$-colinear
map) is a right $A$--linear map $f: M \rightarrow M'$ satisfying $
\rho_{M'} \circ f = (f\tensor{A}\coring{C}) \circ \rho_M$. The
$K$--module of all right $\coring{C}$--colinear maps between two
right comodules $M_{\coring{C}}$ and $M'_{\coring{C}}$ is denoted
by $\homcom{\coring{C}}{M}{M'}$. Right $\coring{C}$--comodules and
their morphisms form a $K$--linear category $\rcomod{\coring{C}}$.
Although not abelian in general, $\rcomod{\coring{C}}$ is a
Grothendieck category provided ${}_A\coring{C}$ is a flat module,
see \cite[Section 1]{ElKaoutit/Gomez/Lobillo:2004c}. The category
$\lcomod{\coring{C}}$ of left $\coring{C}$--comodules is
symmetrically defined.

For more information on corings and comodules, the reader is
referred to \cite{Brzezinski/Wisbauer:2003} and its bibliography.

\section{Exactness of a preradical, localization, and injective objects}\label{localizacion}
In this section we will derive from \cite{Gabriel:1962} some facts
on quotient categories that will be useful in the sequel. Recall
that a full subcategory $\mathcal{C}$ of a Grothendieck category
$\mathcal{G}$ is said to be \emph{closed} if any subobject and any
quotient object of an object belonging to $\cat{C}$ is in
$\cat{C}$, and any direct sum of objects of $\cat{C}$ is in
$\cat{C}$. A closed subcategory $\cat{C}$ of $\cat{G}$ defines a
preradical $\mathfrak{r} : \cat{G} \rightarrow \cat{G}$, which
sends an object $X$ of $\cat{G}$ to its largest subobject
$\mathfrak{r}(X)$ belonging to $\cat{C}$. This preradical is left
exact, since it is right adjoint to the inclusion functor $\cat{C}
\subseteq \cat{G}$. By $\Ker{\mathfrak{r}}$ we denote the full
subcategory of $\cat{G}$ with objects defined by the condition
$\mathfrak{r} (X) = 0$.

A full subcategory $\cat{L}$ of $\cat{G}$ is \emph{dense} if for any
short exact sequence in $\cat{G}$
$$ \xymatrix{0 \ar@{->}[r] & X \ar@{->}[r] & Y \ar@{->}[r] & Z
\ar@{->}[r] & 0 },$$ $Y$ is in $\cat{L}$ if, and only if, both $X$
and $Z$ are in $\cat{L}$. From \cite[15.11]{Faith:1973} we know
that a dense subcategory $\cat{L}$ is $\emph{localizing}$ in the
sense of \cite{Gabriel:1962} if and only if it is stable under
coproducts. Following \cite[Chapter III]{Gabriel:1962}, every
localizing subcategory $\mathcal{L}$ of $\cat{G}$ defines a new
Grothendieck category $\cat{G}/\cat{L}$ (the \emph{quotient}
category), and an exact functor $\mathbf{T}: \mathcal{G}
\rightarrow \mathcal{G}/\mathcal{L}$ (the \emph{canonical
functor}) that admits a right adjoint $\mathbf{S}: \mathcal{G}
/\mathcal{L} \rightarrow \mathcal{G}$. The counit $\phi_{-} :
\mathbf{T} \circ \mathbf{S} \rightarrow
1_{\mathcal{G}/\mathcal{L}} $ of this andjunction is a natural
isomorphism. The unit $\psi_{-}: 1_{\mathcal{G}} \rightarrow
\mathbf{S} \circ \mathbf{T}$ satisfies that both the kernel and
the cokernel of $\psi_X : X \rightarrow (\mathbf{S} \circ
\mathbf{T})(X)$ belong to $\mathcal{L}$ for every object $X$ of
$\mathcal{G}$.

The exactness of a preradical $\mathfrak{r}$ can be expressed in
terms of quotient categories, as the following proposition shows.
The underlying ideas of its proof can be traced back to
\cite[Theorem 2.3]{Gomez/Nastasescu:1995}.

\begin{proposition}\label{nucleo} Let $\cat{C}$ be a closed subcategory of a
Grothendieck category $\cat{G}$ with associated preradical
$\mathfrak{r} : \cat{G} \rightarrow \cat{C}$, and inclusion functor
$\mathfrak{l} : \cat{C} \rightarrow \cat{G}$. The following
statements are equivalent:
\begin{enumerate}[(i)]
\item\label{exacto} $\mathfrak{r}$ is an exact functor;
\item\label{equivalencia} $\cat{K} = \Ker{\fk{r}}$ is a localizing
subcategory of  $\cat{C}$ with canonical functor $\mathbf{T}$, and
there exists an equivalence of categories $H : \cat{G}/\cat{K}
\rightarrow \cat{C}$ such that $H \circ \mathbf{T} =
\mathfrak{r}$.
\end{enumerate}
\end{proposition}
\begin{proof}
(\ref{exacto}) $\Rightarrow$ (\ref{equivalencia}) Since
$\mathfrak{r}$ is exact and preserves coproducts, we easily get
that $\cat{K} = \Ker{\mathfrak{r}}$ is a localizing subcategory.
Consider the canonical adjunctions $$
\xymatrix@R=50pt@C=60pt{\cat{C} \ar@<0,5ex>[r]^{\fk{l}} & \cat{G}
\ar@<0,5ex>[l]^{\fk{r}}}, \qquad \xymatrix@R=50pt@C=60pt{\cat{G}
\ar@<0,5ex>[r]^{\Bf{T}} & \cat{G}/\cat{K} \ar@<0,5ex>[l]^{\Bf{S}}
}$$ where $\Bf{S}$ is right adjoint to $\Bf{T}$, and $\fk{r}$ is
right adjoint to the inclusion functor $\fk{l}$. Composing we get
a new adjoint pair $\Bf{T}\circ \fk{l} : \cat{C} \leftrightarrows
\cat{G}/\cat{K} : \fk{r} \circ \mathbf{S}$, which we claim to
provide an equivalence of categories. The unit of this new
adjunction is given by
$$\xymatrix@C=80pt{ {\rm{id}}_{\cat{C}} = \fk{r} \fk{l} \ar@{->}^-{\fk{r} \,\psi_{\fk{l}}}[r] &
 \fk{r} \, \Bf{S} \, \Bf{T} \, \fk{l} }$$ where $\psi_{-}$ is the
 unit of the adjunction $\Bf{T} \dashv \Bf{S}$. For any object $M$
 of $\cat{C}$, there is an exact sequence
 $$ \xymatrix@C=40pt{ 0 \ar@{->}[r] & X \ar@{->}[r]
 & \fk{l} (M) \ar@{->}^{\psi_{\fk{l}(M)}}[r] & \Bf{S} \Bf{T} \fk{l}(M) \ar@{->}[r] & Y \ar@{->}[r] & 0
 }$$ with $X$ and $Y$ in $\cat{K}$. Apply the exact
 functor $\fk{r}$ to obtain an isomorphism $\fk{r}(\psi_{\fk{l}(M)}):
 M = \fk{r}\fk{l}(M) \,\cong\, \fk{r} \Bf{S} \Bf{T} \fk{l} (M)$.
 Therefore, $\fk{r} \psi_{\fk{l}(-)}$ is a natural isomorphism.
 The counit of the adjunction $\Bf{T}\circ \fk{l} \dashv \fk{r} \circ \mathbf{S}$ is given by the
following composition
 $$\xymatrix@C=40pt{ \Bf{T}\, \fk{l}\, \fk{r} \, \Bf{S} \ar@{->}^{\Bf{T}\, \lambda_{\Bf{S}}}[r]
 & \Bf{T} \, \Bf{S} \ar@{->}_-{\cong}^{\phi_{-}}[r] & {\rm id}_{\cat{G}/\cat{K}}
 }$$ where $\lambda_{-}$ is the counit of the adjunction $\fk{l}
 \dashv \fk{r}$, and $\phi_{-}$ is the counit of the adjunction
 $\Bf{T} \dashv \Bf{S}$. For any object $N$ of $\cat{G}/\cat{K}$,
 $\lambda_{\Bf{S}(N)}$ is a monomorphism with cokernel in
 $\cat{K}$ since $\fk{r}$ is exact. Thus, \cite[Lemme 2, p. 366]{Gabriel:1962} implies that
 $\Bf{T}(\lambda_{\Bf{S}(N)})$ is an isomorphism. Therefore,
 $\phi_N \, \Bf{T}(\lambda_{\Bf{S}(N)})$ is an isomorphism.
 Therefore, $\mathbf{T} \circ \mathfrak{l}$ is an equivalence of
 categories. On the other hand, by \cite[Corollaire 3, p.
 368]{Gabriel:1962}, there exists a functor $H : \cat{G}/\cat{K}
 \rightarrow \cat{C}$ such that $H \circ \mathbf{T} = \mathfrak{r}$.
 By composing on the right with $\mathfrak{l}$ we get $H \circ
 \mathbf{T} \circ \mathfrak{l} = \mathfrak{r} \circ \mathfrak{l} =
 id_{\cat{C}}$. From this, and using that $\mathbf{T} \circ
 \mathfrak{l}$ is an equivalence, we get that $H$ is an equivalence.
 \\
 (\ref{equivalencia}) $\Rightarrow$ (\ref{exacto}) This is obvious,
 since $\mathbf{T}$ is always exact.
\end{proof}

In the rest of this section we consider $\cat{G} = \rmod{B}$, the
category of right modules over a ring $B$.  We fix the following
notation: $\cat{C}$ is a closed subcategory of $\rmod{B}$, with
preradical $\mathfrak{r} : \rmod{B} \rightarrow \cat{C}$, and
inclusion functor $\mathfrak{l} : \cat{C} \rightarrow \rmod{B}$.
We will consider the twosided ideal $\mathfrak{a} =
\mathfrak{r}(B_B)$, and $\cat{K} = \Ker{\mathfrak{r}}$.

\medskip

The following proposition collects a number of well-known
consequences of assuming that $\mathfrak{r}$ is exact. A short
proof is included.

\begin{proposition}\label{plano}
If $\mathfrak{r}$ is exact then $\mathfrak{a}$ is an idempotent
ideal of $B$ such that ${}_B(B/\mathfrak{a})$ is flat, and
$\mathfrak{r}(M) = M \mathfrak{a}$ for every right $B$--module
$M$. In this way, $\cat{K} = \Ker{\mathfrak{r}}$ becomes a
localizing subcategory of $\rmod{B}$ stable under direct products
and injective envelopes.
\end{proposition}
\begin{proof}
Since $\mathfrak{r}$ preserves epimorphisms it follows easily that
$\fk{r}(M) = M \fk{a}$, for any right $B$--module $M$. In
particular, we get that $\cat{K} = \{ M \in \rmod{B} \;|\;
M\mathfrak{a} = 0 \}$. This easily implies that $\cat{K}$ is a
localizing subcategory stable under direct products and essential
extensions. Finally, the flatness of ${}_B(B/\mathfrak{a})$ can be
proved as follows. We know that $\cat{K}$ is isomorphic to
$\rmod{B/\fk{a}}$. Let $\pi: B \rightarrow B/\fk{a}$ be the
canonical projection; the functor $-\tensor{B}(B/\fk{a}) :
\rmod{B} \rightarrow \rmod{B/\fk{a}}$ is left adjoint to the
restriction of scalars functor $\pi_{\ast}: \rmod{B/\fk{a}}
\rightarrow \rmod{B}$. Up to the isomorphism $\cat{K} \cong
\rmod{B/\fk{a}}$, $\pi_{\ast}$ is nothing but the inclusion
functor $j : \cat{K} \rightarrow \rmod{B}$. Since $\cat{K}$ is
stable under injective envelopes, the functor
$-\tensor{B}(B/\fk{a}) $ has to be exact, that is,
${}_B(B/\mathfrak{a})$ is a flat module.
\end{proof}

If $\mathfrak{a}$ is any idempotent ideal of $B$ such that
${}_B(B/\mathfrak{a})$ is flat, then there is a canonical
isomorphism of $B$--bimodules $\mathfrak{a} \cong \mathfrak{a}
\tensor{B} \mathfrak{a}$. This isomorphism makes $\mathfrak{a}$ a
$B$--coring with counit given by the inclusion $\mathfrak{a}
\subseteq B$. We say that $\mathfrak{a}$ is a \emph{left
idempotent $B$--coring} to refer to this situation. The forgetful
functor $U : \rcomod{\mathfrak{a}} \rightarrow \rmod{B}$ induces
then an isomorphism of categories between $\rcomod{\mathfrak{a}}$
and the full subcategory of $\rmod{B}$ whose objects are the
modules $M_B$ such that $M \mathfrak{a} = M$.

\begin{corollary}\label{coanillo}
Assume that $\fk{r}$ is an exact functor. Then
\begin{enumerate}[(i)]
\item The ideal $\fk{a}=\fk{r}(B_B)$ is a left idempotent
$B$--coring whose category of all right comodules
$\rcomod{\fk{a}}$ is isomorphic to the quotient category
$\rmod{B}/\cat{K}$. In particular $\fk{a}$ is a generator of
$\rmod{B}/\cat{K}$. \item The functor $F=\hom{B}{\fk{a}}{-} \circ
\fk{l}: \cat{C} \rightarrow \rmod{B}$ is right adjoint to
$\fk{r}$, where $\fk{l}: \cat{C} \rightarrow \rmod{B}$ is the
inclusion functor. In particular if $E$ is an injective object of
$\cat{C}$, then $F(E)_B$ is an injective right module.
\end{enumerate}
\end{corollary}
\begin{proof}
$(i)$ By Proposition \ref{plano}, $\fk{a}$ is an idempotent $B$-coring. Its category of right comodules clearly
coincides with the torsion class $\cat{C}$, and the stated
isomorphism of categories follows by Proposition \ref{nucleo}. \\
$(ii)$ Given any object $(M,M')$ in $\rmod{B} \times \cat{C}$, we
get natural isomorphisms
\begin{eqnarray*}
\hom{\cat{C}}{\fk{r}(M)}{M'} & \cong & \hom{B}{M\tensor{B}\fk{a}}{\fk{l}(M')} \\
& \cong & \hom{B}{M}{\hom{B}{\fk{a}}{\fk{l}(M')}},
\end{eqnarray*} since $\fk{r}(M) =M \fk{a} \cong M \tensor{B}
\fk{a}$. This means that $F$ is right adjoint to $\fk{r}$. In particular, $F$ preserves
injectives since $\fk{r}$ is exact.
\end{proof}

Given a module $M$ in $\rmod{B}$, the \emph{Wisbauer category}
$\sigma [M]$ associated to $M$ is the full subcategory of
$\rmod{B}$ whose objects are all $M$--subgenerated modules (see
\cite{Wisbauer:1991}). By definition, it is a closed subcategory
and, in fact, it is easy to prove that every closed subcategory of
$\rmod{B}$ is of the form $\sigma [M]$. Therefore, the following
theorem, that summarizes some of the previous results, complements
\cite[42.16]{Brzezinski/Wisbauer:2003}.

\begin{theorem}\label{TeoremaI}
Let $\cat{C}$ be a closed subcategory of a category of modules
$\rmod{B}$ with associated preradical $\fk{r}: \rmod{B} \to
\cat{C}$. Let $\fk{l}: \cat{C} \to \rmod{B}$ be the inclusion
functor, and $\fk{a}=\fk{r}(B)$. The following statements are
equivalent.
\begin{enumerate}[(i)]
\item $\fk{r}: \rmod{B} \rightarrow \cat{C}$ is an exact functor;
\item $\cat{K}=\Sf{Ker}(\fk{r})$ is a localizing subcategory of $\rmod{B}$ with
canonical functor $\mathbf{T}$, and there exists an equivalence $H
: \rmod{B}/\cat{K} \rightarrow \cat{C}$ such that $\mathfrak{r} =
H \circ \mathbf{T}$; \item $F\,=\, \hom{B}{\fk{a}_B}{-} \circ
\fk{l} : \cat{C} \rightarrow \rmod{B}$ is right adjoint to
$\fk{r}$; \item $\fk{a}$ is an idempotent $B$--coring and the
forgetful functor $U : \rcomod{\mathfrak{a}} \rightarrow \rmod{B}$
induces an isomorphism of categories $\rcomod{\coring{a}} \cong
\cat{C}$; \item $\mathfrak{a}^2 = \mathfrak{a}$ and $M\mathfrak{a}
= M$ for every $M$ in $\cat{C}$.
\end{enumerate}
\end{theorem}
\begin{proof}
The equivalences $(i) \Leftrightarrow (ii)$ and $(i) \Leftrightarrow (iii)$ are immediate
from Proposition \ref{nucleo} and Corollary \ref{coanillo}. \\
$(i) \Rightarrow (iv)$ is a consequence of Proposition \ref{nucleo} and Corollary \ref{coanillo}$(i)$.\\
$(iv) \Rightarrow (v)$ Obvious. \\
 $(v) \Rightarrow (i)$ We have easily that $\fk{r}(M)=M\fk{a}$, for every right $B$--module $M$. From
this we get immediately that $\fk{r}$ is a right exact functor.
\end{proof}

Given a right $B$--module $M \in \cat{C}$, by $E_{\cat{C}}(M)$ we
denote its injective hull in the Grothendieck category $\cat{C}$.
According to Theorem \ref{TeoremaI}, if $\mathfrak{r}$ is exact,
then it becomes essentially the canonical functor associated to a
localization with a section functor (the terminology is taken from
\cite{Gabriel:1962}). As a section functor,
$\hom{B}{\mathfrak{a}}{-}$ will preserve injective envelopes, as
stated in Proposition \ref{A-injectivoI}. We give a detailed proof
of this fact, suitable for the forthcoming applications to more
concrete situations.

\begin{proposition}\label{A-injectivoI}
Assume that $\fk{r} : \rmod{B} \rightarrow \cat{C}$ is exact, and
let $M \in \cat{C}$. The map
\begin{equation*}
\zeta_M : M \rightarrow \hom{B}{\mathfrak{a}}{E_{\cat{C}}(M)}
\qquad ( m \mapsto \zeta_M(m)(a) = ma, \; m\in M, a \in
\mathfrak{a})
\end{equation*}
gives an injective envelope of $M$ in $\rmod{B}$. As a
consequence, $M$ is injective in $\rmod{B}$ if and only if $M$ is
injective in $\cat{C}$ and $\zeta_M$ is an isomorphism.
\end{proposition}
\begin{proof}
By Theorem \ref{TeoremaI}, the functor $F =
\hom{B}{\mathfrak{a}}{-} \circ \mathfrak{l} : \cat{C} \rightarrow
\rmod{B}$ is right adjoint to the exact functor $\mathfrak{r}$.
Therefore, $F(E_{\cat{C}}(M)) =
\hom{B}{\mathfrak{a}}{E_{\cat{C}}(M)}$ is injective in $\rmod{B}$.
On the other hand, $\zeta_M$ is obviously a right $B$-linear map.
Let us show that it is injective. Let $m \in M$ such that
$\zeta_M(m) = 0$, that is, $m\mathfrak{a} = 0$. By Theorem
\ref{TeoremaI}, we have $mB=m\mathfrak{a}$, which implies $m = 0$.
Let us prove that $\zeta_M$ is essential. Pick a non zero element
$f \in \hom{B}{\mathfrak{a}}{E_{\cat{C}}(M)}$, so there exists $0
\neq u \in \fk{a}$ such that $0 \neq f(u) \in E_{\cat{C}}(M)$.
Since $M$ is essential in $E_{\cat{C}}(M)$, there exists a non
zero element $b \in B$ such that $0 \neq f(u)b \in M$. Since
$B/\mathfrak{a}$ is flat as a left $B$--module (Proposition
\ref{plano}), there exists $ w\in \fk{a}$ with $ub = ubw$ (see,
e.g. \cite[42.5]{Brzezinski/Wisbauer:2003}). If we consider the
map $g=fub$, then $g(x)=(fub)(x)=f(ubx)$, for all $x \in \fk{a}$,
that is $g=fub =\zeta_M(f(ub))$ is a non zero element of
$\zeta_M(M)$, as $g(w) = f(ub) \neq 0$.
\end{proof}

\begin{definition}\label{supfin}
Assume that $\mathfrak{a}$ has a set of local units in the sense
of \cite{Abrams:1983}, that is, $\mathfrak{a}$ contains a set $E$
of commuting idempotents such that for every $x \in \mathfrak{a}$
there exits $e \in E$ such that $xe = ex = x$. A right $B$--module
$M$ is said to be of \emph{finite support} if there exits a finite
subset $F \subseteq E$ such that $e \in E$ and $m\sum_{e \in F}e =
m$ for every $m \in M$.
\end{definition}

A straightforward argument proves that if $M_B$ is of finite
support, then every $f \in \hom{B}{\mathfrak{a}}{M}$ is of the
form $f(x) = mx$ for some $m \in M$. Therefore, we deduce from
Proposition \ref{A-injectivoI}:

\begin{corollary}\label{injCinjA}
Assume that $\mathfrak{a}$ has a set of local units, and let $M
\in \cat{C}$ of finite support. Then $M$ is injective in
$\rmod{B}$ if and only if $M$ is injective in $\cat{C}$. As a
consequence, given a homomorphism of rings $A \rightarrow B$ with
${}_AB$ flat, we deduce that if $M$ is injective in $\cat{C}$,
then $M$ is injective in $\rmod{A}$.
\end{corollary}

\begin{remark}
In Definition \ref{supfin} and Corollary \ref{injCinjA}, it
suffices to assume that $\mathfrak{a}$ contains a set of commuting
idempotents $E$ such that $\mathfrak{a} = \sum_{e \in E}eB$.
\end{remark}

In what follows we specialize our results to the case where the
subcategory $\cat{C}$ is isomorphic to the category of right
comodules over a given $A$-coring $\coring{C}$. This is the case
when $\coring{C}$ is member of a rational pairing
$\TT=(\coring{C},B,\esc{-}{-})$. Rational pairings for coalgebras
over commutative rings were introduced in \cite{Gomez:1998} and
used in \cite{Abuhlail/Gomez/Lobillo:2001} to study the category
of right comodules over the finite dual coalgebra associated to
certain algebras over Noetherian commutative rings. This
development was adapted for corings in
\cite{ElKaoutit/Gomez/Lobillo:2004c}, see also
\cite{Abuhlail:2003a}.

Recall from \cite[Section 2]{ElKaoutit/Gomez/Lobillo:2004c} that a
three-tuple $\TT=(\coring{C},B,\esc{-}{-})$ consisting of an
$A$-coring $\coring{C}$, an $A$-ring $B$ (i.e., $B$ is an algebra
extension of $A$) and a balanced $A$--bilinear form $\esc{-}{-}:
\coring{C} \times B \to A$, is said to be \emph{a right rational
pairing over} $A$ provided
\begin{enumerate}[(1)]
 \item $\beta_A : B \to {}^*\coring{C}$ is a ring anti-homomorphism, where ${}^*\coring{C}$ is
the left dual convolution ring of $\coring{C}$ defined in
\cite[Proposition 3.2]{Sweedler:1975}, and

\item $\alpha_M$ is an injective map, for each right $A$-module $M$,
\end{enumerate}

where $\alpha_{-}$ and $\beta_{-}$ are the following  natural transformations
$$\xymatrix@R=0pt{\beta_N : B \tensor{A} N \ar@{->}[r] &
\hom{}{_{A}\coring{C}}{{}_{A}N}, \\ b \tensor{A} n \ar@{->}[r]& \left[ c
\mapsto \esc{c}{b}n \right] }\quad \xymatrix@R=0pt{ \alpha_M : M
\tensor{A} \coring{C} \ar@{->}[r] & \hom{}{B_{A}}{M_{A}} \\
m \tensor{A} c \ar@{->}[r] & \left[ b \mapsto m\esc{c}{b} \right].
}$$

Given a right rational pairing $\TT=(\coring{C},B,\esc{-}{-})$
over $A$, we can define a functor called the \emph{right rational
functor} as follows.  An element $m$ of a right $B$--module $M$ is
called \emph{rational} if there exists a set of \emph{right
rational parameters} $\{(c_i,m_i)\} \subseteq \coring{C} \times M$
such that $ m b = \sum_i m_i\esc{c_i}{b}$, for all $b \in B$. The
set of all rational elements in $M$ is denoted by $\Rat^{\TT}(M)$.
As it was explained in \cite[Section
2]{ElKaoutit/Gomez/Lobillo:2004c}, the proofs detailed in
\cite[Section 2]{Gomez:1998} can be adapted in a straightforward
way in order to get that $\Rat^{\TT}(M)$ is a $B$--submodule of
$M$ and the assignment $M \mapsto \Rat^{\TT}(M)$ is a well defined
functor
\begin{equation*}
\Rat^{\TT} : \rmod{B} \rightarrow \rmod{B},
\end{equation*}
which is in fact a left exact preradical. Therefore, the full
subcategory $\Rat^{\TT}(\rmod{B})$ of $\rmod{B}$ whose objects are
those $B$--modules $M$ such that $\Rat^{\TT}(M) = M$ is a closed
subcategory. Furthermore, $\Rat^{\TT}(\rmod{B})$ is a Grothendieck
category which is shown to be isomorphic to the category of right
comodules $\rcomod{\coring{C}}$ as \cite[Theorem
2.6']{ElKaoutit/Gomez/Lobillo:2004c} asserts (see also
\cite[Proposition 2.8]{Abuhlail:2003a}).

\begin{example}\label{par-canonico}
Let $\coring{C}$ be an $A$--coring such that ${}_A\coring{C}$ is a
locally projective left module (see \cite[Theorem
2.1]{Zimmermann-Huisgen:1976} and \cite[Lemma
1.29]{Abuhlail:2003a}). Consider the endomorphism ring
$\lend{\coring{C}}{\coring{C}}$ as a subring of the endomorphism
ring $\lend{A}{\coring{C}}$, that is, with multiplication opposite
to the composition of maps. Since $\Delta_{\coring{C}}$ is a left
$\coring{C}$--colinear and a right $A$--linear map, the canonical
ring extension $A \rightarrow \lend{A}{\coring{C}}$ factors
throughout the extension $\lend{\coring{C}}{\coring{C}}
\hookrightarrow \lend{A}{\coring{C}}$. Therefore, the three-tuple
$\TT=(\coring{C},\lend{\coring{C}}{\coring{C}},\esc{-}{-})$, where
the balanced $A$-bilinear $\esc{-}{-}$ map is defined by
$\esc{c}{f}= \varepsilon_{\coring{C}}(f(c))$, for $(c,f) \in
\coring{C} \times \lend{\coring{C}}{\coring{C}}$ is a rational
pairing since $\lend{\coring{C}}{\coring{C}}$ is already a ring
anti-isomorphic to ${}^*\coring{C}$ via the beta map associated to
$\esc{-}{-}$. We refer to $\TT$ as \emph{the right canonical
pairing} associated to $\coring{C}$.
\end{example}

The following theorem complements
\cite[20.8]{Brzezinski/Wisbauer:2003}.

\begin{theorem}\label{Teorema}
Let $\TT=(\coring{C},B,\esc{-}{-})$ be a right rational pairing
with rational functor $\Rat^{\TT}: \rmod{B} \rightarrow \rmod{B}$,
and put $\mathfrak{a}= \Rat^{\TT}(B_B)$,
$\cat{K}=\Ker{\Rat^{\TT}}$. The following statements are
equivalent:
\begin{enumerate}[(i)]
\item $\Rat^{\TT}: \rmod{B} \rightarrow \rmod{B}$ is an exact
functor; \item $\cat{K}$ is a localizing subcategory of $\rmod{B}$
with canonical functor $\mathbf{T}$, and there exists an
equivalence $H : \rmod{B}/\cat{K} \rightarrow
\Rat^{\TT}(\rmod{B})$ such that $\Rat^{\TT} = H \circ \mathbf{T}$;
\item $F\,=\, \hom{B}{\fk{a}_B}{-} \circ \fk{l} :
\Rat^{\TT}(\rmod{B}) \rightarrow \rmod{B}$ is right adjoint to
$\Rat^{\TT}$; \item $\fk{a}$ is an idempotent $B$--coring and the
forgetful functor $U : \rcomod{\mathfrak{a}} \rightarrow \rmod{B}$
induces an isomorphism of categories $\rcomod{\fk{a}} \cong
\Rat^{\TT}(\rmod{B}) \cong \rcomod{\coring{C}}$; \item
${}_B\fk{a}$ is a pure submodule of ${}_BB$, $\fk{a}^2 = \fk{a}$,
and $\coring{C} \fk{a} = \coring{C}$.
\end{enumerate}
\end{theorem}
\begin{proof}
By Theorem \ref{TeoremaI} we only need to show that $(v)
\Rightarrow (iv)$ since $(iv) \Rightarrow (v)$ is clear. We have
that $\fk{a}$ is an idempotent $B$--coring and $\coring{C} \cong
\coring{C} \tensor{B} \fk{a}$ as right $B$--modules. Given any
rational right $B$--module $X$ with its canonical structure of
right $\coring{C}$--comodule, we obtain a $B$--linear isomorphism
$X \cong X\cotensor{\coring{C}}\left( \coring{C} \tensor{B}
\fk{a}\right)$ (recall that the comultiplication is a right
$B$--linear map), where the symbol $-\cotensor{\coring{C}}-$
refers to the cotensor bifunctor over $\coring{C}$. Using the left
version of \cite[Lemma 2.2]{Gomez:2002}, we get
$$X\,\cong\, X\cotensor{\coring{C}} \coring{C} \, \cong\,
X\cotensor{\coring{C}}\left( \coring{C} \tensor{B} \fk{a}\right)
\,\cong \, \left( X\cotensor{\coring{C}}\coring{C}
\right)\tensor{B} \fk{a} \, \cong \, X\tensor{B} \fk{a}.$$  That
is, $X$ is in fact a right $\fk{a}$--comodule.
\end{proof}

\begin{remark}
Right rational pairings are instances of right coring measurings
in the sense of \cite{Brzezinski:2004a}. In this way, given an
exact rational functor $\Rat^{\TT}$ the isomorphism of categories
$\rcomod{\coring{C}} \cong \rcomod{\fk{a}}$ stated in Theorem
\ref{Teorema} can be interpreted as an isomorphism of corings in
an adequate category. Following to \cite[Definition
2.1]{Brzezinski:2004a}, a $B$--coring $\coring{D}$ is called a
\emph{right extension} of an $A$--coring $\coring{C}$ provided
$\coring{C}$ is a $(\coring{C},\coring{D})$--bicomodule with the
left regular coaction $\Delta_{\coring{C}}$. Corings understood as
pairs $(\coring{C}:A)$ (i.e., $\coring{C}$ is an $A$--coring) and
morphisms understood as right coring extensions (i.e., a pairs
consisting of an action and coaction) with their bullet
composition form a category denoted by $\mathbf{CrgExt}^r_K$ (see
\cite{Brzezinski:2004a} for more details). If we apply this to the
setting of Corollary \ref{Teorema}, then it can be easily checked
that $(\coring{C}: A)$ and $(\fk{a}:B)$ become isomorphic objects
in the category $\mathbf{CrgExt}^r_K$.
\end{remark}

From Proposition \ref{A-injectivoI} and Corollary \ref{injCinjA},
we obtain:

\begin{proposition}\label{A-injectivo}
Let $\TT=(\coring{C},B,\esc{-}{-})$ be a right rational pairing
with rational functor $\Rat^{\TT}: \rmod{B} \rightarrow \rmod{B}$,
and put $\mathfrak{a}= \Rat^{\TT}(B_B)$. Assume that $\Rat^{\TT}$
is an exact functor. Let $M$ be a right $\coring{C}$--comodule,
and $E(M_{\coring{C}})$ its injective hull in
$\rcomod{\coring{C}}$.
\begin{enumerate}[(a)]
\item The map
\begin{equation*}
\zeta_M : M \rightarrow \hom{B}{\mathfrak{a}}{E(M_{\coring{C}})}
\qquad ( m \mapsto \zeta_M(m)(a) = ma, \; m\in M, a \in
\mathfrak{a})
\end{equation*}
gives an injective envelope of $M$ in $\rmod{B}$. %
\item $M$ is injective in $\rmod{B}$ if and only if $M$ is
injective in $\rcomod{\coring{C}}$ and $\zeta_M$ is an
isomorphism. %
\item If ${}_AB$ is flat, $\zeta_M$ is an isomorphism, and $M$ is
injective in $\rcomod{\coring{C}}$, then $M$ in injective in
$\rmod{A}$. %
\item If $\mathfrak{a}$ has a set of local units and $M$ is of
finite support, then $M$ is injective in $\rcomod{\coring{C}}$ if
and only if $M$ in injective in $\rmod{B}$. %
\item Assume that $\mathfrak{a}$ has a set of local units, $M$ is
of finite support and ${}_AB$ is flat. If $M$ is injective in
$\rcomod{\coring{C}}$, then $M$ is injective in $\rmod{A}$.
\end{enumerate}
\end{proposition}

\section{Rational functors for entwined and Doi-Koppinen modules}

In this section, we shall study the exactness of the rational
functors for the corings coming from entwining structures. When
particularized to the entwining structures given by a comodule
algebra, we will obtain a result from \cite{Abuhlail:2003a}. Most
of results in \cite{Dascalescu/Nastasescu/Torrecillas:1999} are
deduced.

\subsection{Entwining structures with rational
functor}\label{ent-funtRat}

Recall from \cite{Brzezinski/Majid:1998} that an entwining
structure over $K$ is a three-tuple $(A,C)_{\psi}$ consisting of a
$K$--algebra $A$ with multiplication $\mu$ and unity $1$, a
$K$--coalgebra $C$ with comultiplication $\Delta$ and counity
$\varepsilon$, and a $K$--module map $\psi: C\tensor{K}A
\rightarrow A\tensor{K}C$ satisfying
\begin{equation}\label{est-entw}
\begin{array}{l}
\psi \circ (\tensork{C}{\mu}) = (\tensork{\mu}{C}) \circ
(\tensork{A}{\psi}) \circ (\tensork{\psi}{A}), \\ \\
(\tensork{A}{\Delta}) \circ \psi = (\tensork{\psi}{C}) \circ
(\tensork{C}{\psi}) \circ (\tensork{\Delta}{A}), \\ \\ \psi \circ
(\tensork{C}{1}) = \tensork{1}{C}, \hspace{3.em}
(\tensork{A}{\varepsilon}) \circ \psi = \tensork{\varepsilon}{A}.
\end{array}
\end{equation}
By \cite[Proposition 2.2]{Brzezinski:2002} the corresponding
$A$--coring is $\coring{C} = A\tensor{K} C$ with the $A$--bimodule
structure given by $a''(a'\tensor{K}c)a=a''a'\psi(c\tensor{K}a)$,
$a,a',a'' \in A$, $c \in C$, the comultiplication
$\Delta_{\coring{C}}= A\tensor{K}\Delta$, and the counit
$\varepsilon_{\coring{C}}= A\tensor{K}\varepsilon$. Furthermore,
the category of right $\coring{C}$--comodules is isomorphic to the
category of right entwined modules.

The map $(\phi,\nu): (C,K) \rightarrow (\coring{C},A)$ defined by
$\nu(1)=1$ and $\phi(c)=1\tensor{K}c$, is a homomorphism of
corings in the sense of \cite{Gomez:2002}. As in \cite{Gomez:2002}
the associated induction and ad-induction functors to this
morphism are, respectively, given by $\cat{O}: \rcomod{\coring{C}}
\rightarrow \rcomod{C}$ and $-\tensor{K}A: \rcomod{C} \rightarrow
\rcomod{\coring{C}}$, where $\cat{O}$ is the cotensor functor
$-\cotensor{\coring{C}}(A\tensor{K}C)$. When $\rcomod{\coring{C}}$
is interpreted as the category of entwined modules, $\cat{O}$ is
naturally isomorphic to the forgetful functor.  Moreover, there is
a natural isomorphism
$$\xymatrix@R=0pt{ \hom{\coring{C}}{M \tensor{K}A}{N}
\ar@{->}^-{\cong}[r]& \hom{C}{M}{\cat{O}(N)}
\\ f \ar@{|->}[r] & [m \mapsto f(m\tensor{K}1)] \\
[m \tensor{K} a \mapsto g(m)a] & g, \ar@{|->}[l] }$$ for every
pair of comodules $(M_C,\, N_{\coring{C}})$. Thus the functor
$\cat{O}$ is a right adjoint functor of $-\tensor{K}A$. If $C_K$
is a flat module, then $\cat{O}$ is exact, since $U_A:
\rcomod{\coring{C}} \rightarrow \rmod{A}$ is already an exact
functor (see, \cite[Proposition
1.2]{ElKaoutit/Gomez/Lobillo:2004c}).

We know from \cite{Brzezinski:2002} that the left dual convolution
ring ${}^*\coring{C}$ is isomorphic as a $K$--module to
$\hom{K}{C}{A}$. Up to this isomorphism the convolution
multiplication reads
\begin{equation}\label{pro-smach}
f \cdot g \,=\, \mu \circ (A\tensor{K}f)\, \circ\, \psi\, \circ\,
(C \tensor{K} g) \, \circ \, \Delta_C, \quad f, g \in
\hom{K}{C}{A}.
\end{equation} The connection between this convolution ring and the
usual coalgebra convolution ring $C^*$ is given by the following
homomorphism of rings
\begin{equation}\label{phiestr}
\xymatrix@R=0pt{\Phi: C^* \ar@{->} [r] & {}^*\coring{C}, & & ( x
\ar@{|->}[r] & \tensork{A}{x} ).}
\end{equation}

\begin{proposition}\label{lratentwing}
Let $(A,C)_{\psi}$ be an entwining structure over $K$ such that
$C_K$ is a locally projective module and consider its
corresponding $A$--coring $\coring{C}=A\tensor{K}C$. Suppose that
there is a right rational pairing $\T=(\coring{C},B,\esc{-}{-})$
and an anti-morphism of $K$--algebras $\varphi: C^* \rightarrow B$
which satisfy the following two conditions: (1) $\beta \, \circ \,
\varphi = \Phi$, where $\beta: B \rightarrow {}^*\coring{C}$ is
the anti-homomorphism of $K$--algebras associated to $\T$ and
$\Phi$ is the homomorphism of rings given in equation
\eqref{phiestr}; (2) for every pair of elements $(a,x) \in A\times
C^*$, there exists a finite subset of pairs $\{(x_i,a_i)\}_i
\subseteq C^* \times A$ such that $a\varphi(x) \, = \, \sum_i
\varphi(x_i) a_i$. Then, by restricting scalars we have
$$\Rat^{\T}\left(M_B\right) \, = \,
\Rat^r_C\left({}_{C^*}M\right)$$ for every right $B$--module $M$,
where $\Rat^r_{C}(-)$ is the canonical right rational functor
associated to the $K$--coalgebra $C$.
\end{proposition}
\begin{proof}
Start with an arbitrary element $m \in \Rat^{\T}(M_B)$ with right
rational system of parameters
$\{(\sum_k\tensork{a_{kj}}{c_{kj}},\,m_j)\}_j \subset \coring{C}
\times M$. Then for every $x \in C^*$, we have
\begin{eqnarray*}
  xm  \,=\, m \varphi(x) &=& \sum_{k,j} m_j
  \esc{a_{kj} \tensor{K}c_{kj}}{\varphi(x)} \\
   &=& \sum_{k,j} m_j
\beta(\varphi(x))(a_{kj} \tensor{K} c_{kj}) \\
   &=& \sum_{k,j}
m_j \Phi(x)(a_{kj} \tensor{K} c_{kj}),\quad \Phi = \beta \circ
\varphi, \\
   &=&  \sum_{k,j}m_j a_{kj} x(c_{kj}),
\end{eqnarray*} thus $\{c_{kj},m_ja_{kj}\} \subset C \times M $ is
a right rational system of parameters for $m \in {}_{C^*}M$; that
is  $m \in \Rat^r_C({}_{C^*}M)$. Therefore, $\Rat^{\T}(M_B)
\subseteq \Rat^r_C({}_{C^*}M)$. Conversely, start with a pair of
elements $(a,x) \in A\times C^*$, and let $\{(x_i,a_i)\}_i \subset
C^* \times A$ be the finite system given by hypothesis, that is
$a\varphi(x) = \sum_i \varphi(x_i)a_i$. So, for every element $m
\in \Rat^{r}_C({}_{C^*}M)$ with right $C$--coaction
$\rho_{\Rat^{r}_C({}_{C^*}M)}(m) = \sum_{(m)} m_{(0)} \tensor{K}
m_{(1)}$, we have
\begin{eqnarray*}
x(ma) \,=\, (ma) \varphi(x) &=& m(a\varphi(x)) \\ &=& \sum m
(\varphi(x_i) a_i),\,\, a \varphi(x) = \sum \varphi(x_i) a_i
\\ &=& \sum (x_i m) a_i
\\ &=& \sum (m_{(0)} x_i(m_{(1)}))a_i \\ &=& \sum m_{(0)}
\left(\underset{}{}\Phi(x_i)(1 \tensor{K} m_{(1)})a_i \right)\\
&=& \sum m_{(0)} \left(\underset{}{} \beta(\varphi(x_i))(1 \tensor{K} m_{(1)})a_i\right) \\
&=& \sum m_{(0)} \left(\underset{}{} \beta (\varphi(x_i)a_i)(
1\tensor{K} m_{(1)})\right), \,\, \beta \text{ is right }
A-\text{linear} \\ &=&
\sum m_{(0)} \, \esc{1 \tensor{K} m_{(1)}}{\varphi(x_i)a_i} \\
&=& \sum m_{(0)}\, \esc{1 \tensor{K}m_{(1)}}{a\varphi(x)} \\ &=&
\sum m_{(0)}\, \esc{a_{\psi} \tensor{K}
m_{(1)}^{\psi}}{\varphi(x)} ,\,\, \psi(m_{(1)} \tensor{K} a)= \sum
a_{\psi} \tensor{K} m_{(1)}^{\psi} \\ &=& \sum m_{(0)} a_{\psi}
x(m_{(1)}^{\psi}).
\end{eqnarray*}
We conclude that $ma \in \Rat^r_C({}_{C^*}M)$ with right
$C$--coaction $\rho_{\Rat^r_C({}_{C^*}M)}(ma) = \sum
m_{(0)}a_{\psi} \tensor{K} m_{(1)}^{\psi}$. From which we conclude
that $\Rat^r_C({}_{C^*}M)$ is an entwined module, and thus a right
$\coring{C}$--comodule or, equivalently, a right rational
$B$--submodule of $M_B$. Therefore, $\Rat^r_C({}_{C^*}M) \subseteq
\Rat^{\T}(M_B)$.
\end{proof}

\subsection{The category of Doi-Koppinen modules}

We apply the results of the subsection \ref{ent-funtRat} to the
category of Doi-Koppinen modules. This category is identified with
the category of right rational modules over a well known ring.
Some results of this section were proved by different methods for
particular Hopf algebras in \cite[Theorem 2.3]{Cai/Chen:1994}, for
algebras and coalgebras over a field in \cite[Proposition
2.7]{Dascalescu/Nastasescu/Torrecillas:1999}, and more recently
for bialgebras in \cite[Theorem 3.18, Proposition
3.21]{Abuhlail:2003a}.

\bigskip

Let $\h$ be a Hopf $K$--algebra, $(A,\rho_A)$ a right
$\h$--comodule $K$--algebra, and $(C,\varrho_C)$ a left
$\h$--module $K$--coalgebra. That is $\rho_A: A \rightarrow
A\tensor{K}\h$ and $\varrho_C: \h\tensor{K}C \rightarrow C$ are,
respectively, a $K$--algebra and a $K$--coalgebra map. We will use
Sweedler's notation, that is $\Delta_C(c) =\sum_{(c)}
c_{(1)}\tensor{K}c_{(2)}$,
$\Delta_{\h}(h)=\sum_{(h)}h_{(1)}\tensor{K}h_{(2)}$, and
$\rho_A(a)=\sum_{(a)}a_{(1)}\tensor{K}a_{(2)}$, for every $c \in
C$, $h \in\h$ and $a \in A$.

Following \cite{Doi:1992, Koppinen:1995}, a \emph{Doi-Koppinen
module} is a left $A$--module $M$ with a structure of right
$C$--comodule $\rho_M$ such that, for every $a \in A$, $m \in M$,
$$\rho_M (am) = \sum a_{(0)}m_{(0)} \tensor{K} a_{(1)}m_{(1)}.$$
A morphism between two Doi-Koppinen modules is a left $A$--linear
and right $C$--colinear map. Doi-Koppinen modules and their
morphisms form the category $\dkmod{A}{C}$.

Consider the following $K$--map ($A^{o}$ means the opposite ring
of $A$)
\begin{equation}\label{dk-psi}
\xymatrix@R=0pt{ \psi: C \tensor{K} A^{o} \ar@{->}[r] & A^{o}
\tensor{K} C, & (c \tensor{K} a^{o} \ar@{|->}[r] & \sum_{(a)}
a_{(0)}^o \tensor{K} a_{(1)}c).}
\end{equation}
It is easily seen that the map $\psi$ satisfies all identities of
equation \eqref{est-entw}. That is $(A^o,C)_{\psi}$ is an
entwining structure over $K$. So, consider the associated
$A^o$--coring $\coring{C} = A^o\tensor{K}C$, the $A^o$--biactions
are then given by $$ b^o(a^o \tensor{K} c)\, =\, (ab)^o \tensor{K}
c,\, \text{ and }\, (b^o\tensor{K} c)a^o \,=\, b^o \psi(c
\tensor{K} a^o)=\sum (a_{(0)}b)^o \tensor{K} a_{(1)}c,$$ for every
$a^o, b^o \in A^o$ and $c \in C$. The convolution multiplication
of $\hom{K}{C}{A^o}$ comes out from the general equation
\eqref{pro-smach}, as
\begin{equation}\label{producto}
f.g(c)\,=\,\sum_{(c)} \left( f(g(c_{(2)})_{(1)}c_{(1)})\,
g(c_{(2)})_{(0)} \right)^o \,\, \in \, A^o,\, f,g \in
\hom{K}{C}{A^{o}} \text{ and } c \in C.
\end{equation} This multiplication coincides with the generalized smash
product of $A$ by $C$, denoted by $\sharp(C,A)$ in
\cite[(2.1)]{Koppinen:1995}.

Define the smash product $A\sharp C^*$ whose underling $K$--module
is the tensor product $A\tensor{K}C^*$ and internal multiplication
is given by $$(a \sharp x).(b \sharp y)\,=\,\sum ab_{(0)} \sharp
(xb_{(1)})y,$$ for $a\tensor{K}x,\,b\tensor{K}y \in A\tensor{K}
C^*$, and where the left $\h$--action on $C^*$ is induced by the
right $\h$--action on $C$. The unit of this multiplication is $1
\sharp \varepsilon_{C}$. Moreover, its clear that the maps
$$\xymatrix@R=0pt@C=40pt{ -\sharp \varepsilon_C:\,A \ar@{->}[r] & A\sharp C^*,\\
a \ar@{|->}[r] & a\sharp \varepsilon_C } \quad
\xymatrix@R=0pt@C=40pt{ 1\sharp -:\,C^* \ar@{->}[r] & A\sharp C^*
\\ x \ar@{|->}[r] & 1 \sharp x }$$ are $K$--algebra maps, and an
easy computation shows that $$\xymatrix@R=0pt{ \alpha_A: A\sharp
C^* \ar@{->}[r] & \hom{K}{C}{A^o}, & (a \sharp x\ar@{|->}[r] & [c
\mapsto a^ox(c)]) }$$ is also a $K$--algebra morphism where
$\hom{K}{C}{A^o}$ is endowed with the multiplication of equation
\eqref{producto}.

\begin{proposition}\cite[Proposition
3.21]{Abuhlail:2003a}\label{dkpar} Let $\h$ be a Hopf
$K$--algebra, $A$ a right $\h$--comodule $K$--algebra and $C$ a
left $\h$--module $K$--coalgebra. Consider
$\coring{C}=A^o\tensor{K}C$ the $A^o$--coring associated to the
entwining structure $(A^o,C)_{\psi}$ where $\psi$ is defined by
\eqref{dk-psi}, and let $B = (A \sharp C^*)^o$. Suppose that $C_K$
is a locally projective module. Then
$\T=(\coring{C},B,\esc{-}{-})$ is right rational pairing over
$A^o$ with the bilinear form $\esc{-}{-}$ defined by
$$\xymatrix@R=0pt{\coring{C} \times B \ar@{->}[r] & A^{o} \\ (a^o
\tensor{K}c,( b \sharp x)^o) \ar@{|->}[r] & \esc{a^o
\tensor{K}c}{(b \sharp x)^o}\,=\,a^ob^ox(c) }$$ $a,b\in A$, $c \in
C$, $x \in C^*$. Moreover, $(A^o,C)_{\psi}$, $\T$ and
$\varphi=(1\sharp -)^o: C^* \rightarrow B$ satisfy the conditions
(1) and (2) stated in Proposition \ref{lratentwing}, and using
restriction of scalars, we obtain
$$ \Rat^{\T}(M_B) \,=\, \Rat^r_C({}_{C^*}M),$$ for every right
$B$--module $M$.
\end{proposition}
\begin{proof}
First we show that $\esc{-}{-}$ is bilinear and balanced. For $a,
b,e\in A$, $x \in C^*$ and $c \in C$, we compute
\begin{eqnarray*}
  \esc{a^o\tensor{K}c}{(b\sharp x)^o e^o} &=& \esc{a^o\tensor{K}c}{((e\sharp \varepsilon_C) (b\sharp x))^o} \\
   &=& \sum \esc{a^o\tensor{K}c}{(eb_{(0)}\sharp (\varepsilon_Cb_{(1)})x)^o}  \\
   &=&\sum a^o\, b_{(0)}^o\, e^o\, ((\varepsilon_C\,b_{(1)})x)(c) \\
   &=&\sum a^o\, b_{(0)}^o\, e^o\, \varepsilon_{\h}(b_{(1)})x(c) \\
   &=& a^o\, b^o\, e^o \, x(c) \\
   &=& \esc{a^o\tensor{K}c}{(b\sharp x)^o }\, e^o,
\end{eqnarray*} which shows that $\esc{-}{-}$ is right
$A^o$-linear, and
\begin{eqnarray*}
  \esc{(a^o\tensor{K}c)e^o}{(b\sharp x)^o } &=& \esc{a^o\psi(c\tensor{K}e^o)}{(b\sharp x)^o} \\
   &=&\sum \esc{a^oe_{(0)}^o\tensor{K}e_{(1)}c}{(b\sharp x)^o} \\
   &=&\sum a^o e_{(0)}^o b^o x(e_{(1)}c)  \\
   &=&\sum a^o e_{(0)}^o b^o \,\, (xe_{(1)})(c)  \\
   &=& \esc{a^o\tensor{K}c}{ \sum(be_{(0)})\sharp xe_{(1)}}  \\
   &=& \esc{a^o\tensor{K}c}{(b\sharp x)(e\sharp\varepsilon_C)} \\
   &=& \esc{a^o\tensor{K}c}{e^o\,(b\sharp x)^o },
\end{eqnarray*} which proves that $\esc{-}{-}$ is $A^o$-balanced.
The pairing $\esc{-}{-}$ is clearly left $A^o$--linear. Consider
now the right natural transformation associated to $\esc{-}{-}$:
$$\xymatrix@R=0pt{\alpha_N : \hspace{0.5em} N
\tensor{A^{o}}\coring{C} \ar@{->}[r] & \hom{A^{o}}{B_{A^{o}}}{N}
\\ n \tensor{A^{o}} a^o \tensor{K} c \ar@{|->}[r]& \left[(b \sharp x)^o
\mapsto n \esc{a^o \tensor{K} c}{(b \sharp
x)^o}\,=\,n(a^ob^ox(c))\right]. }
$$ We need to show that $\alpha_{-}$ is injective. So let $\sum_i
n_i \tensor{A^o} 1 \tensor{K} c_i \in \, N\tensor{A^o}\coring{C}$
whose image by $\alpha_{N}$ is zero. Since $C_K$ is locally
projective, associated to the finite set $\{c_i\}_i$ there exists
a finite set $\{(c_l,\,x_l)\} \subset C \times C^*$ such that $c_i
= \sum_l c_lx_l(c_i)$. The condition $$\alpha_N(\sum_i n_i
\tensor{A^{o}}1 \tensor{K}c_i)((1 \sharp x_{l})^o) \,\,=\,\,
\sum_i n_i x_{l}(c_i) =0, \text{ for all the }l's, $$ implies that
$$\sum_i n_i \tensor{A^{o}}1 \tensor{K}c_i = \sum_{i,\,l}
n_i\tensor{A^{o}}1 \tensor{K} x_{l}(c_i)c_{l} = \sum_{l}
\left(\sum_i n_i x_{l}(c_i)\right) \tensor{A^{o}} 1 \tensor{K}
c_{l} =0,$$ That is $\alpha_N$ is an injective map for every right
$A^o$--module $N$. Therefore, $\T$ is a right rational system.
Lastly, the map $\beta: B \rightarrow {}^*\coring{C}$ sending
$(b\sharp x)^o \mapsto [a^o\tensor{K}c \mapsto
\esc{a^o\tensor{K}c}{(b\sharp x)^o}]$ is an anti-homomorphism of
$K$--algebras, and $B$ is a $K$--algebra extension of $A^o$, thus
$\T$ is actually a right rational pairing.

Let $a^o \in A^o$, $c \in C$ and $x \in C^*$, then
$$\beta(\varphi(x))(a^o\tensor{K}c)=\esc{a^o\tensor{K}c}{(1\sharp
x)^o} \,=\, a^ox(c) \,=\,\Phi(x)(a^o\tensor{K}c)$$ which implies
the condition $(1)$ of Proposition \ref{lratentwing}. For the
condition $(2)$, it is easily seen that the set
$\{a_{(0)}^o,\,xa_{(1)}\}$, where $\rho_A(a)= \sum a_{(0)}
\tensor{K}a_{(1)}$, satisfies this condition for the pair $(a^o,x)
\in A^o \times C^*$. The last stated assertion is a consequence of
Proposition \ref{lratentwing}, and this finishes the proof.
\end{proof}

\begin{theorem}\label{ultimo}
Let $\h$ be a Hopf $K$--algebra, $A$ a right $\h$--comodule
$K$--algebra and $C$ a left $\h$--module $K$--coalgebra. Consider
$\coring{C}=A^o\tensor{K}C$ the $A^o$--coring associated to the
entwining structure $(A^o,C)_{\psi}$ where $\psi$ is defined by
\eqref{dk-psi}, and set $B = (A \sharp C^*)^o$. Suppose that $C_K$
is a locally projective module and consider the right rational
pairing $\T=(\coring{C},B,\esc{-}{-})$ over $A^o$ of Proposition
\ref{dkpar}, and put $\fk{a}=\Rat^{\T}(B_B)$. If $A_K$ is a flat
module and $\Rat^r_C(-)$ is an exact functor, then
\begin{enumerate}[(a)]
\item $\Rat^r_C({}_{C^*}C^*)$ is a right $\h$--submodule of
$C^*$ and $\fk{a}\,= A\tensor{K}\Rat^r_C({}_{C^*}C^*)$. %
\item For each right $B$--module $M$, the map
\begin{equation}\label{isorat}
\xymatrix{\hom{B}{\mathfrak{a}_B}{M} \ar[r] &
\hom{C^*}{\Rat^r_C({}_{C^*}C^*)}{M}}
\end{equation}
sending $f$ onto the morphism $\widehat{f}$ defined by
$\widehat{f}(c^*) = f(1 \tensor{} c^*)$ for $c^* \in C^*$ is an
isomorphism of $K$--modules. %
\item If, for every left $A \sharp C^*$--module $M$, we endow
$\hom{C^*}{\Rat_C^r ({}_{C^*}C^*)}{M}$ with the structure of a
left $A \sharp C^*$--module transferred from that of $\hom{A
\sharp C^*}{{}_{A \sharp C^*}\mathfrak{a}}{M}$ via the isomorphism
(\ref{isorat}), then we obtain a functor
\begin{equation*}
\xymatrix{\hom{C^*}{\Rat_C^r ({}_{C^*}C^*)}{-} : \dkmod{A}{C}
\ar[r] & \lmod{A \sharp C^*} }
\end{equation*}
which is right adjoint to the functor (see Proposition \ref{dkpar})
\begin{equation*}
\xymatrix{\Rat^r_C :  \lmod{A \sharp C^*} \ar[r] & \dkmod{A}{C}}
\end{equation*}
\end{enumerate}
\end{theorem}
\begin{proof}
$(a)$ Let $y \in \Rat^r_C({}_{C^*}C^*)$ with rational system of
parameters  $\{(y_i,c_i)\}_i \subset C^* \times C$. For any $h \in
\h$ and $x \in C^*$, we obtain as in \cite[Lemma
3.1]{Dascalescu/Nastasescu/Torrecillas:1999}:
\begin{eqnarray*}
  (x(yh))(c) &=& \sum_{(c)} x(c_{(1)}) y(hc_{(2)}) \\
   &=& \sum_{(c),(h)} x(\varepsilon_{\h}(h_{(1)})c_{(1)}) y(h_{(2)}c_{(2)}) \\
   &=& \sum_{(c),(h)} x(S(h_{(1)})h_{(2)}c_{(1)}) y(h_{(3)}c_{(2)}) \\
   &=& \sum_{(h)} (xS(h_{(1)})y)(h_{(2)}c) \\
   &=& \sum_{(h),i} y_i(h_{(2)}c) (xS(h_{(1)}))(c_i) \\ &=& \sum_{(h),i} y_i(h_{(2)}c) x(S(h_{(1)})c_i),
\end{eqnarray*} for every $c \in C$, where $S$ is the antipode of
$\h$. That is, $x(yh) = \sum_{(h),\,
i}(y_ih_{(2)})\,x(S(h_{(1)})c_i)$. Hence,
$\{(y_ih_{(2)},S(h_{(1)})c_i)\} \subset C^* \times C$ is a
rational system of parameters for $yh$. Thus $yh \in
\Rat^r_C({}_{C^*}C^*)$, and $\Rat^r_C({}_{C^*}C^*)$ is a right
$\h$--submodule of $C^*$. Since $A_K$ is a flat module, an easy
computation shows now that $A\tensor{K}\Rat^r_C({}_{C^*}C^*)$ is a
two-sided ideal of $A\sharp C^*$. Let $x \in C^*$, $a\tensor{K}y
\in A\tensor{K}\Rat^r_C({}_{C^*}C^*)$, and $\{(y_i,c_i)\}_i
\subset C^* \times C$ a rational system of parameters for $y$.
Applying the smash product, we get $$x (a \tensor{K}
y)\,=\,\sum_{(a)} a_{(0)} \tensor{K}((xa_{(1)})y)
\,=\,\sum_{(a),\,i} (a_{(0)} \tensor{K} y_i) x(a_{(1)}c_i);$$ this
means that $\{(a_{(0)}\tensor{K}y_i, a_{(1)}c_i)\}_{(a),\,i}
\subset (A\tensor{K}\Rat^r_C({}_{C^*}C^*))\times C$ is a rational
system of parameters for $a\tensor{K}y \in
{}_{C^*}\lr{A\tensor{K}\Rat^r_C({}_{C^*}C^*)}$. Proposition
\ref{dkpar}, implies now that $A\tensor{K}\Rat^r_C({}_{C^*}C^*)
\subseteq \fk{a}$. Conversely, we know that $\fk{a}$ is a right
$\coring{C}$--comodule, so the underlying $K$-module is a right
$C$--comodule, and, since $\Rat^r_C$ is exact, $\fk{a} =
\Rat^r_C({}_{C^*}C^*) \fk{a}$. From this equality, it is easy to
see that $\fk{a} \subseteq A \tensor{K}\Rat^r_C({}_{C^*}C^*)$, and
the desired equality is
derived.\\
\noindent $(b)$ We know that $B = (A \sharp C^*)^{o
}$ and, by
$(a)$, we have $\mathfrak{a} = A \tensor{K} \Rat^r({}_{C^*}C^*)$.
Consider the homomorphism of algebras $C^* \rightarrow A \sharp
C^*$, which gives, as usual, the induction functor $(A \sharp C^*)
\tensor{C^*} - : \lmod{C^*} \rightarrow \lmod{A \sharp C^*}$ which
is left adjoint to the restriction of scalars functor $\lmod{A
\sharp C^*} \rightarrow \lmod{C^*}$. The mapping $f \mapsto
\widehat{f}$ is then defined as the composition
\begin{multline*}
\hom{A \sharp C^*}{A \tensor{K} \Rat^r_C({}_{C^*}C^*)}{M} \cong \\
\hom{A \sharp C^*}{(A\sharp C^*)\tensor{C^*}
\Rat^r_C({}_{C^*}C^*)}{M}  \cong
\hom{C^*}{\Rat^r_C({}_{C^*}C^*)}{M},
\end{multline*}
where the second is the adjointness isomorphism, and the first one
comes from the obvious isomorphism $(A \sharp C^*) \tensor{C^*}
\Rat^r_C({}_{C^*}C^*) \cong A \tensor{K} \Rat^r_C({}_{C^*}C^*)$.\\
\noindent $(c)$ This is a consequence of $(b)$ and Theorem
\ref{Teorema}.
\end{proof}

Keep, in the following corollary, the hypotheses of Theorem
\ref{ultimo}.

\begin{corollary}
If $M$ is an object of $\dkmod{A}{C}$, and $E({}_AM^C)$ denotes
its injective hull in the category $\dkmod{A}{C}$, then
\begin{enumerate}[(a)]
\item The map
\begin{equation*}
\zeta_M : M \rightarrow \hom{C^*}{
\Rat^r_C({}_{C^*}C^*)}{E({}_AM^C)} \qquad ( m \mapsto
\zeta_M(m)(c^*) = c^*m, \; m\in M))
\end{equation*}
gives an injective envelope of $M$ in $\lmod{A \sharp C^*}$. %
\item $M$ is injective in $\lmod{A \sharp C^*}$ if and only if $M$
is injective in $\dkmod{A}{C}$ and $\zeta_M$ is an
isomorphism. %
\item Assume that the antipode of $\h$ is bijective, and that
$C^*$ is flat as a $K$--module. Let $M$ be injective in
$\dkmod{A}{C}$. If $M$ has finite support as a right
$C$--comodule, then $M$ is injective as a left $A$--module.
\end{enumerate}
\end{corollary}
\begin{proof}
The two first statements follow from Proposition \ref{A-injectivo}
and Theorem \ref{ultimo}. For the last statement, observe that
$B_{A^{op}}$ is a flat module. Now, the proof of \cite[Lemma
2.6]{Dascalescu/Nastasescu/Torrecillas:1999} runs here to prove
that ${}_{A^{o}}B$ is flat.
\end{proof}

\begin{remark}
We have proved that, under suitable conditions,
\begin{equation}\label{t}
\Rat^{\T}(M_B)=\Rat^r_C({}_{C^*}M)\,=\,M \fk{a} \,=\,
\Rat^r_C({}_{C^*}C^*)M,
\end{equation}
for every right $B$--module $M$. Therefore, equation \eqref{t}
establishes a radical functor: $t: \lmod{A\sharp C^*} \rightarrow
\dkmod{A}{C}$ which acts on objects by $M \rightarrow
\Rat_C^r({}_{C^*}C^*)M$. This radical was used in \cite[Lemma
2.9]{Dascalescu/Nastasescu/Torrecillas:1999} for left and right
semiperfect coalgebras over a commutative field. In this way, if we
apply our results and \cite[Proposition 2.2]{Gomez/Nastasescu:1995}
to this setting, then most part of the results stated in
\cite{Dascalescu/Nastasescu/Torrecillas:1999} become consequences of
the results stated in this paper. In particular, let us mention
that, for a semiperfect coalgebra over a field, any comodule of
finite support in the sense of
\cite{Dascalescu/Nastasescu/Torrecillas:1999} becomes of finite
support in the sense of Definition \ref{supfin}. Finally, let us
note that the results from
\cite{Dascalescu/Nastasescu/Torrecillas:1999} are only applicable to
group-graded algebras over a field. This restriction has been
dropped by our approach, and we fully cover the case of graded rings
(take $K=\mathbb{Z}$), since the $\mathbb{Z}$--coalgebra
$\mathbb{Z}G$, where the elements of the group $G$ are all
group-like, is easily shown to have an exact rational functor. Of
course, the category of comodules over this coalgebra is not
semiperfect.
\end{remark}

\providecommand{\bysame}{\leavevmode\hbox
to3em{\hrulefill}\thinspace}
\providecommand{\MR}{\relax\ifhmode\unskip\space\fi MR }
\providecommand{\MRhref}[2]{%
  \href{http://www.ams.org/mathscinet-getitem?mr=#1}{#2}
} \providecommand{\href}[2]{#2}


\end{document}